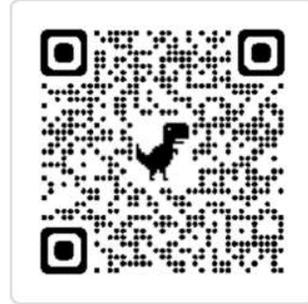

# Two-Person Stochastic Duel with Energy Fuel Constraint Ammo

SONG-KYOO (AMANG) KIM

**ABSTRACT**

This paper deals with a novel variation of the versatile stochastic duel game that incorporates an energy fuel constraint into a two-player duel game. The energy fuel not only measures the vitality of players but also determines the power of the shooting projectile. The game requires players to carefully balance their energy usage, while trying to outmaneuver their opponent. This unique theoretical framework for the stochastic game model provides a valuable method for understanding strategic behavior in competitive environments, particularly in decision-making scenarios with fluctuating processes. The proposed game provides players with the challenge of optimizing their energy fuel usage, while managing the risk of losing the game. This novel model has potential for implementation across diverse fields, as it allows for a versatile conception of energy fuel. These energy fuels may encompass conventional forms, such as natural gas, petroleum, and electrical power, and even financial budgets, human capital, and temporal resources. The unique rules and constraints of the game in this research are expected to contribute insights into the decision-making strategies and behaviors of players in a wide range of practical applications. This research primarily focuses on deriving compact closed-form solutions, utilizing transformation and flexible analysis techniques adapted to varying the concept of the energy fuel level. By presenting a comprehensive description of our novel analytical approach and its application to the proposed model, this study aims to elucidate the fundamental principles underlying the energy fuel constraint stochastic duel game model.

# 1. INTRODUCTION

Game theory is a powerful tool for modeling and analyzing strategic decision-making in a wide range of fields, from economics and political science to computer science and engineering [1-6]. Modern game theory has seen significant developments in recent years, particularly in the areas of repeated games, evolutionary game theory, and network games. One important recent development is the use of machine learning techniques in game theory, which has led to advances in predicting and optimizing outcomes in complex, multi-player games [7]. Another area of active research is the analysis of games in which players have incomplete information, which has applications in contract theory, mechanism design, and more [8]. Network games have also become an important area of research, as they model interactions among agents in a social or economic network, with implications for contagion and diffusion dynamics [9]. Another growing area of interest is the use of game theory in analyzing cybersecurity, particularly in the context of defending against cyber-attacks [3, 5, 10]. Overall, modern game theory has seen exciting developments and applications in a wide range of fields. As researchers continue to develop new methods and models, game theory is likely to remain a valuable tool for understanding strategic decision-making in complex, dynamic environments. Duel games have been the focus of study in game theory due to their applicability in modeling competitive situations, ranging from military conflicts to economic decisions [11]. Recent research has explored variations of the duel game, such as the generalized stochastic duel game [12] and the duel game with multiple asymmetric players [13, 14] and penalties. Other studies have examined the impact of different strategies on the outcome of duel games, including the use of mixed strategies [15] and time dependent strategies [16, 17]. Several studies have also explored the application of duel games in specific contexts, such as decision-making in the presence of incomplete information [18], analysis of political competitions [19], and modeling of predator-prey interactions in ecological systems [20].

# 2. STOCHASTIC DUEL GAME WITH FUEL CONSTRAINT AMMO

The antagonistic duel game of two players (called "A" and "B") in the time domain is considered and both players know the full information regarding the success probabilities based on the time domain [1, 2]. Let us assigned that $P_a(t)$ is the monotone increasing CPFs (Probability Density Functions) for player A regarding hitting the opponent player (player B) at the time $t$. Similarly, $P_b(t)$ is the probability of player B for hitting the opponent player (player A) at the time $t$. Each hitting probability (CDF) could be arbitrary chosen and it reaches 1 (i.e., 100 %) when the time $t$ goes to the infinite ($< \infty$). The strategic decision for the shooting moment means finding the moment when a player will have the best chance to hit the other. There is a certain point that maximizes the chance for succeeding in the shoot, and this optimal point becomes the moment of success in the continuous time domain. This moment $t^*$ is defined as follows:

$$t^* = inf\{t \geq 0 : P_a(t) + P_b(t) \geq 1\}. \qquad (2.1)$$

Additionally, each player has the energy fuel level which randomly drooped and the spent (or usage) of the energy fuel for each player is monotone and non-decreasing. Let ($\Omega$,

$\mathcal{F}(\Omega),P)$ be probability space $\mathcal{F}_A$, $\mathcal{F}_B$, $\mathcal{F}_\tau \subseteq \mathcal{F}(\Omega)$ be independent $\sigma$-subalgebras. Suppose:

$$\mathcal{A} := \sum_{k \geq 0} \varepsilon_{s_k}, \ s_0(=0) < s_1 < s_2 < \cdots, \ \text{a.s.} \tag{2.2}$$

$$\mathcal{B} := \sum_{j \geq 0} \varepsilon_{t_j}, \ t_0(=0) < t_1 < t_2 < \cdots, \ \text{a.s.} \tag{2.3}$$

are $\mathcal{F}_A$-measurable and $\mathcal{F}_B$-measurable renewal point processes ($\varepsilon_a$ is a point mass at $a$) with intensities $\lambda_a$ and $\lambda_b$ and position and point independent marking. These two values are related with the energy drain of each player. The energy level of each player becomes the shooting power to hit a opponent player. The energy fuel of player A drains at times $s_1, s_2, \ldots$ and the magnitudes of energy drains are formalized by process $\mathcal{A}$. The energy drain of player B are described by process $\mathcal{B}$ similarly. The processes $\mathcal{A}$ and $\mathcal{B}$ are specified their transforms:

$$\mathbb{E}\left[g^{\mathcal{A}(s)}\right] = e^{\lambda_a(s)(g-1)}, \ \mathbb{E}\left[h^{\mathcal{B}(t)}\right] = e^{\lambda_b(t)(h-1)}. \tag{2.4}$$

The game is observed at random times in accordance with the point process

$$\mathcal{T} := \sum_{i \geq 0} \varepsilon_{\tau_i}, \ \tau_0(>0)), \tau_1, \ldots, \tag{2.5}$$

which is assumed to be delayed renewal process.

$$(A(t), B(t)) := \mathcal{A} \otimes \mathcal{B}([0, \tau_k]), \ k = 0, 1, \ldots, \tag{2.6}$$

forms an observation process upon $\mathcal{A} \otimes \mathcal{B}$ embedded over $\mathcal{T}$, with respective increments

$$(X_k, Y_k) := \mathcal{A} \otimes \mathcal{B}([\tau_{k-1}, \tau_k]), \ k = 1, 2, \ldots, \tag{2.7}$$

and

$$X_0 = A_0, \ Y_0 = B_0. \tag{2.8}$$

The observation process could be formalized as

$$\mathcal{A}_\tau \otimes \mathcal{B}_\tau := \sum_{k \geq 0} (X_k, Y_k) \varepsilon_{\tau_k}, \tag{2.9}$$

where

$$\mathcal{A}_\tau = \sum_{i \geq 0} X_i \varepsilon_{\tau_i}, \ \mathcal{B}_\tau = \sum_{i \geq 0} Y_i \varepsilon_{\tau_i}, \tag{2.10}$$

and it is with position dependent marking and with $X_k$ and $Y_k$ being dependent with the notation

$$\Delta_k := \tau_k - \tau_{k-1}, \ k = 0, 1, \ldots, \ \tau_{-1} = 0, \tag{2.11}$$

and

$$\sigma(z,\theta) = \mathbb{E}\left[z^{(X_i-Y_i)^+}e^{-\theta\Delta_i}\right], z > 0. \tag{2.12}$$

By using the double expectation, we have

$$\sigma(z,\theta) = \sigma(\theta + (\lambda_a - \lambda_b)(1-z)) \tag{2.13}$$

where

$$\sigma(\theta) = \mathbb{E}\left[e^{-\theta\Delta_i}\right], \sigma_0(\theta) = \mathbb{E}\left[e^{-\theta\tau_0}\right]. \tag{2.14}$$

Let us consider the maximum energy levels of players $M_a$ and $M_b$. The energy level of player A after draining fuel is $M_a - \mathcal{A}_\tau$ and $M_b - \mathcal{B}_\tau$ for player B from (2.9). The stochastic process the energy level of each player are as follows:

$$\{M_a - \mathcal{A}_\tau\} \otimes \{M_b - \mathcal{B}_\tau\} := (M_a, M_b) - \sum_{k\geq 0}(X_k, Y_k)\varepsilon_{\tau_k}, \tag{2.15}$$

and the game is over when the $k$-th observation epoch $\tau_k$, the shooting power of player A which is equivalent with the remained energy level at the moment of the shooting is greater than the energy fuel level of player B:

$$M_a - \mathcal{A}_{\tau_k} \geq M_b - \mathcal{B}_{\tau_k}. \tag{2.16}$$

To further formalization of the game, the exit index could be defined as follows:

$$\nu := \inf\{k : A_k = A_0 + X_1 + \cdots + X_k \leq \Sigma_a\}, \tag{2.17}$$

$$\mu := \inf\{j : B_j = H_0 + Y_1 + \cdots + Y_j \leq \Sigma_b\}. \tag{2.18}$$

where

$$\Sigma_a = M_a - (M_b - B_j),$$

$$\Sigma_b = M_b - (M_a - A_k).$$

Since player A is assumed to win the game at time $\tau_\nu$, we shall be targeting the confined game in the view point of player A. The passage tme $\tau_\nu$ is associated exit time from the confined game and the formula (2.15) will be modified as

$$\left\{M_a - \overline{\mathcal{A}_\tau}\right\} \otimes \left\{M_b - \overline{\mathcal{B}_\tau}\right\} := (M_a, M_b) - \sum_{k\geq 0}^{\nu}(X_k, Y_k)\varepsilon_{\tau_k}, \tag{2.19}$$

which the path of the game from $\mathcal{F}(\Omega) \cap \{A_\nu - B_\nu \leq |M_a - M_b|\} \cap \{\mathcal{T} \geq t^*\}$, which gives an exact definition of the model observed until $\tau_\nu$. The joint functional of the stochastic duel with the fuel limited ammunition is as follows:

$$\Phi_{\{\nu,\mu\}} = \Phi_{\{\nu(\Sigma_a),\mu(\Sigma_b)\}}(\zeta, z_0, z_1, \theta_0, \theta_1) \tag{2.20}$$

$$= \mathbb{E}\left[\zeta^{\nu} z_0^{(A_{\nu-1}-B_{\mu-1})} z_1^{(A_{\nu}-B_{\mu})} e^{-\theta_0 \tau_{\nu-1}} e^{-\theta_1 \tau_{\nu}} \mathbf{1}_{\{\nu-\mu \leq \frac{\mathcal{M}_{ab}}{\sigma}\}} \mathbf{1}_{\{\nu \geq \frac{t^*}{\sigma}\}}\right]$$

where

$$\{A_{\nu} - B_{\mu} \leq |M_a - M_b|\} \cap \{\mathcal{T} \geq t^*\} \equiv \left\{\nu - \mu \leq \frac{\mathcal{M}_{ab}}{\sigma}\right\} \cap \left\{\nu \geq \frac{t^*}{\sigma}\right\}, \quad (2.21)$$

$$\mathcal{M}_{ab} = |M_a - M_b|, \overline{\sigma} = \mathbb{E}[\Delta_k]. \quad (2.22)$$

This model will represent the status of both players upon the *exit time* $\tau_{\nu}$ and the pre-exit time $\tau_{\nu-1}$. The pre-exit time is a particular interest because player A wants to predict not only her time for the highest chance, but also the moment for the next highest chance prior to this. The **Theorem-1** establishes an explicit formula for $\Phi_{\{\mu,\nu\}}$ and we abbreviate with (2.23)-(2.27):

$$\gamma(z,\theta) = \sigma(z,\theta), \gamma_0(z,\theta) = \sigma_0(z,\theta), \quad (2.23)$$

$$\phi_A(x,\theta) = \mathbb{E}\left[x^{X_k} e^{-\theta \Delta_j}\right] = \sigma(\theta - \lambda_a(1-x)), \quad (2.24)$$

$$\phi_A^0(x,\theta) = \mathbb{E}\left[x^{A_0} e^{-\theta \tau_0}\right] = \sigma_0(\theta - \lambda_a(1-x)), \quad (2.25)$$

$$\phi_B(y,\theta) = \mathbb{E}\left[y^{Y_j} e^{-\theta \Delta_j}\right] = \sigma(\theta - \lambda_b(1-y)), \quad (2.26)$$

$$\phi_B^0(y,\theta) = \mathbb{E}\left[y^{B_0} e^{-\theta \tau_0}\right] = \sigma_0(\theta - \lambda_b(1-y)). \quad (2.27)$$

The linear operators are defined as follows:

$$\mathcal{D}_{(p,q)}[f(p,q)](x,y) := (1-x)(1-y) \sum_{p \geq 0} \sum_{q \geq 0} f(p,q) x^p y^q, \quad (2.28)$$

then

$$f(p,q) = \mathfrak{D}_{(x,y)}^{(p,q)}\left[\mathcal{D}_{(p,q)}\{f(p,q)\}\right], \quad (2.29)$$

where $\{f(p,q)\}$ is a sequence, with the inverse

$$\mathfrak{D}_{(x,y)}^{(p,q)}(\bullet) = \begin{cases} \left(\frac{1}{p! \cdot q!}\right) \lim_{(x,y) \to 0} \frac{\partial^p \partial^q}{\partial x^p \partial y^q} \frac{1}{(1-x)(1-y)}(\bullet), & p \geq 0, q \geq 0, \\ 0, & \text{otherwise.} \end{cases} \quad (2.30)$$

**Theorem-1.** The functional $\Phi_{\mu\nu}$ of the game on trace $\sigma$-algebra
$\mathcal{F}(\Omega) \cap \{A_{\nu} - B_{\nu} \leq |M_a - M_b|\} \cap \{\mathcal{T} \geq t^*\}$ satisfies the following formula:

$$\Phi_{\{\nu,\mu\}} = \mathfrak{D}_{(x,y)}^{(M_a, M_b)} \left\{ \frac{\zeta \Gamma_0 \gamma_1 (1 - \phi_x)(\zeta \Gamma)^{\left\lceil \frac{t^*}{\sigma} \right\rceil}}{\phi_y^{\left\lfloor \frac{M_{ab}}{\sigma} \right\rfloor} (1 - \zeta \Gamma)} \right\}, \qquad (2.31)$$

where
$$\Gamma = \gamma(z_0 z_1 x, \theta_0 + \theta_1), \qquad (2.32)$$

$$\Gamma_0 = \gamma_0(z_0 z_1 x, \theta_0 + \theta_1), \qquad (2.33)$$

$$\gamma_1 = \gamma(z_1, \theta_1), \qquad (2.34)$$

$$\phi_x = \phi_A(x, 0), \phi_y = \phi_B(y, 0), \qquad (2.35)$$

***Proof:*** We find the explicit formula of the joint function $\Phi_{\{\mu,\nu\}}$ which starts from (2.19):

$$\Phi_{\{\nu,\mu\}} = \Phi_{\{\nu(\Sigma_a),\mu(\Sigma_b)\}}(\zeta, z_0, z_1, \theta_0, \theta_1)$$

$$= \sum_{j \geq 0} \sum_{k \geq 0} \zeta^j \mathbb{E}\left[ z_0^{(A_{\nu-1} - B_{\nu-1})} z_1^{(A_\nu - B_\nu)} e^{-\theta_0 \tau_{\nu-1}} e^{-\theta_1 \tau_\nu} \right.$$

$$\left. \mathbf{1}_{\left\{ j-k \leq \left\lfloor \frac{M_{ab}}{\sigma} \right\rfloor \right\}} \mathbf{1}_{\left\{ j \geq \left\lceil \frac{t^*}{\sigma} \right\rceil \right\}} \mathbf{1}_{\{\nu=j, \mu=k\}} \right]$$

and applying the operator $\mathcal{D}$ to random family $\{\mathbf{1}_{\{\nu(x)=p, \mu(y)=q\}} : x \geq 0\}$, we have

$$\mathcal{D}_{(p,q)}\left[\mathbf{1}_{\{\nu(x)=k, \mu(y)=j\}}\right](x, y) = \left(x^{A_{k-1}} - x^{A_k}\right)\left(y^{B_{j-1}} - y^{B_j}\right) \qquad (2.36)$$

and from the previous research [1-5],

$$\Psi(x, y) = \sum_{j \geq 0} \sum_{k \geq 0} \zeta^j \mathcal{D}_{(x,y)} \mathbb{E}\left[ z_0^{(A_{\nu-1} - B_{\nu-1})} z_1^{(A_\nu - B_\nu)} e^{-\theta_0 \tau_{\nu-1}} e^{-\theta_1 \tau_\nu} \right. \qquad (2.37)$$

$$\left. \mathbf{1}_{\left\{ j-k \leq \left\lfloor \frac{M_{ab}}{\sigma} \right\rfloor \right\}} \mathbf{1}_{\left\{ j \geq \left\lceil \frac{t^*}{\sigma} \right\rceil \right\}} \mathbf{1}_{\{\nu=j, \mu=k\}} \right]$$

$$= \sum_{j \geq \left\lceil \frac{t^*}{\sigma} \right\rceil} L_{1j} L_{2j} \sum_{k \geq j - \left\lfloor \frac{M_{ab}}{\sigma} \right\rfloor} L_{3jk} L_{4jk}. \qquad (2.38)$$

where
$$L_{1j} = \zeta^j \mathbb{E}\left[ (z_0 z_1)^{(A_{j-1} - B_{j-1})} (xy)^{B_{j-1}} e^{-(\theta_0 + \theta_1)\tau_{j-1}} \right] = \zeta^j \Gamma_0 \Gamma^{j-1}, \qquad (2.39)$$

$$L_{2j} = \mathbb{E}\left[ z_1^{(X_j - Y_j)} e^{-\theta_1 \Delta_j} (1 - x^{X_j}) \right] = \gamma_1 (1 - \phi_x), \qquad (2.40)$$

$$L_{3jk} = \mathbb{E}[y^{Y_j+Y_{j+1}+\cdots+Y_{k-1}}] = \phi_y^{k-j}, \tag{2.41}$$

$$L_{4jk} = \mathbb{E}[1 - y^{Y_k}] = 1 - \phi_y. \tag{2.42}$$

From (2.39)-(2.42),

$$\Psi(x,y) = \sum_{j \geq \lceil \frac{t^*}{\sigma} \rceil} (\zeta^j \Gamma_0 \Gamma^{j-1})(\gamma_1(1 - \phi_x)) \left( \phi_y^{-\lceil \frac{\mathcal{M}_{ab}}{\sigma} \rceil} \right)$$

$$= \left( \frac{\zeta \Gamma_0 \gamma_1 (1 - \phi_x)}{\phi_y^{\lceil \frac{\mathcal{M}_{ab}}{\sigma} \rceil}} \right) \left( \frac{(\zeta \Gamma)^{\lceil \frac{t^*}{\sigma} \rceil}}{1 - \zeta \Gamma \phi_{xy}} \right). \tag{2.43}$$

Therefore,

$$\Psi(x,y) = \frac{\zeta \Gamma_0 \gamma_1 (1 - \phi_x)(\zeta \Gamma)^{\lceil \frac{t^*}{\sigma} \rceil}}{\phi_y^{\lceil \frac{\mathcal{M}_{ab}}{\sigma} \rceil}(1 - \zeta \Gamma)}. \tag{2.44}$$

From (2.29) and (2.44), finally we have:

$$\Phi_{\{\nu,\mu\}} = \mathfrak{D}_{(x,y)}^{(M_a, M_b)} \left\{ \frac{\zeta \Gamma_0 \gamma_1 (1 - \phi_x)(\zeta \Gamma)^{\lceil \frac{t^*}{\sigma} \rceil}}{\phi_y^{\lceil \frac{\mathcal{M}_{ab}}{\sigma} \rceil}(1 - \zeta \Gamma)} \right\}. \tag{2.36}$$

The functional $\Phi_{\{\nu,\mu\}}$ contains all decision making parameters regarding this standard stopping game. The information includes the optimal number of iterations of players (i.e., $\nu$ and $\mu$), the best moments of shooting ($\tau_\nu$, $\tau_\mu$; *exit time*) and the one step prior to the best shooting times ($\tau_{\nu-1}$, $\tau_{\mu-1}$; *pre-exit time*). The information for player A from the closed functional are as follows:

$$\mathbb{E}[\nu] = \lim_{\zeta \to 1} \left( \frac{\partial}{\partial \zeta} \right) \Phi_{\{\nu,\mu\}}(\zeta, 1, 1, 0, 0), \tag{2.37}$$

$$\mathbb{E}[\tau_{\nu-1}] = \lim_{\theta \to 0} \left( -\frac{\partial}{\partial \theta} \right) \Phi_{\{\nu,\mu\}}(1, 1, 1, \theta, 0). \tag{2.38}$$

Additionally, there are some special case could be considered. First case is the assets of both players are the same (i.e., $\mathcal{M}_{ab} = 0$). From (2.36), the formula is changed as follows:

$$\Phi_{\{\mu,\nu\}} = \mathfrak{D}_{(x,y)}^{(M_a,M_b)} \left\{ \frac{\zeta \Gamma_0 \gamma_1 (1-\phi_x)(\zeta \Gamma)^{\left\lceil \frac{t^*}{\sigma} \right\rceil}}{1-\zeta\Gamma} \right\}, \mathcal{M}_{ab} = 0. \qquad (2.39)$$

This implication means that both players have the same energy fuel level and the winning strategy for this game should take the shot as soon as passing the threshold from (2.1). The other case might be the initial energy level of player A is smaller than what player B has (i.e., $M_a < M_b$). The best strategy of player A should wait until player B shoots (and fail) because player A has no chance to win even he hits the target correctly.

## 3. Conclusion

A new successor to the antagonistic stochastic duel game was studied. This research primarily focused on deriving compact closed-form solutions utilizing transformation and flexible analysis techniques, which were adapted by varying the concept of the energy fuel level. In this innovative duel game, a player can win the game only if their bullets hit the target player and if their shooting power exceeds the remaining energy fuel level of the target player at the moment of shooting. A joint functional of the standard stopping game was constructed to analyze the strategic decision parameters, which indicated the best moment for shooting in the time domain stochastic game. This study provides a thorough account of the innovative analytic approach, shedding light on the fundamental principles that underlie the stochastic duel game model with energy fuel constraints.